\documentclass[11pt]{amsart}
\usepackage{amsmath}
\usepackage{amssymb}
\usepackage{enumerate}

\newtheorem{thm}{Theorem}
\newtheorem{lem}[thm]{Lemma}
\newtheorem{prop}[thm]{Proposition}
\newtheorem{cor}[thm]{Corollary}
\newtheorem{rmk}{Remark}

\begin{document}

\title{Recurrence times and large deviations}

\author{Yong Moo Chung}

\address{Department  of  Applied  Mathematics, 
Hiroshima  University,
Higashi-Hiroshima 739-8527, Japan}

\curraddr{Department of Mathematics,
The Pennsylvania State University \\
University Park, PA 16802, USA}

\email{chung@amath.hiroshima-u.ac.jp}

\date{November 21, 2006. {\it Last modified:} December 14, 2007.\\
\hspace{0.4cm}{\it 2000 Mathematics Subject Classification.} 
Primary 37D25; Secondary 37D35, 60F10.}

\keywords{large deviations, reccurence times, nonuniform hyperbolicity.}

\thanks{This research was supported by
Grant-in-Aid for Young Scientists (B) of JSPS, Grant No.17740059}

\begin{abstract}
We give a criterion to determine the large deviation rate functions
for abstract dynamical systems on towers.
As an application of this criterion  we show the level 2
large deviation principle for some class of smooth interval maps 
with nonuniform hyperbolicity.
\end{abstract}

\pagestyle{plain} 

\maketitle

\section{Introduction}

Let $I$ be a compact metric space with a finite Borel measure $m$
as a reference measure.
Unless otherwise stated, $m$ will be normalized Lebegue measure 
if $I$ is a manifold.
We denote by $\mathcal{M}$ the space of the Borel probability measures
on $I$ equipped with the weak* topology.
For a nonsingular transformation $f:I\to I$, not necessary invariant
for $m$, 
we say that it satisfies {\it the (level 2) large deviation principle}
if  
there is an upper semicontinuous function
$q: \mathcal{M} \to [-\infty , 0] ,$
called \textit{the rate function,}
satisfying
$$\liminf_{n\to\infty}
\frac{1}{n}\log m\left(\left\{ x\in I :
\delta_x^n \in \mathcal{G} \right\}\right)
\ge \sup_{\mu\in \mathcal{G}} q(\mu )$$
for each open set 
$\mathcal{G}\subset\mathcal{M} ,$
and
$$\limsup_{n\to\infty}
\frac{1}{n}\log m\left(\left\{ x\in I : 
\delta_x^n  \in \mathcal{C} \right\}\right)
\le \max_{\mu\in \mathcal{C}} q(\mu )$$
for each closed set
$\mathcal{C}\subset\mathcal{M} ,$
respectively, where 
$$\delta_x^n := \frac{1}{n} \sum_{i=0}^{n-1} \delta_{f^i (x)} 
\in\mathcal{M}$$
denotes the empirical distribution along the orbit of $f$ through $x\in I$.
We refer to Ellis' book \cite{E} 
for a general theory of large deviations and 
its background in statistical mechanics.

It is well-known that for a uniformly 
or partially hyperbolic dynamical system on a manifold
with a specification property
it satisfies the large deviation principle
and the rate function is represented
as the difference between the metric entropy and 
the sum of positive Lyapunov exponents \cite{OP,Y1}.
A similar result as above is also known for piecewise expanding maps, 
and then the rate function $q$ coincides with 
the free energy function $F$ given by
$$F(\mu ) :=
\begin{cases}
 \displaystyle h_{\mu} (f) - \int\log \vert f' \vert d\mu ,
 & \text{ for } \mu \in \mathcal{M}_f ,
 \\
\displaystyle  -\infty & \text{ otherwise,}
\end{cases}$$
where $\mathcal{M}_f$ denotes the set of 
$f$-invariant Borel probability measures,
and $h_{\mu} (f)$ the metric entropy of $\mu\in\mathcal{M}_f$
for $f$ \cite{T1}.
The results above on the large deviation principle
include Ruelle's inequality, 
Pesin's and Rohlin's formulas for entropy
\cite{Ke, Pe, Ru}.

Some of the large deviations estimates are also known 
for nonuniformly hyperbolic dynamical systems.
Keller and Nowicki \cite{KN}
gave a large deviations theorem 
for a nonrenormalizable unimodal map $f:I\to I$ 
satisfying the Collet-Eckmann condition 
that:
for any continuous function $\varphi$ of bounded variation
with positive variance
\begin{equation}
\label{ldp:KN}
\alpha (\varepsilon ) :=\lim_{n\to\infty}
\frac{1}{n}\log m \left(\left\{ 
x\in I : | \frac{1}{n}S_n\varphi (x) -
\int\varphi d\mu_0  |\ge \varepsilon 
\right\}\right) <0
\end{equation}
exists for small $\varepsilon >0 ,$
where $$S_n \varphi (x):= \varphi (x)+ \varphi (f(x))+\cdots 
+\varphi (f^{n-1}(x)),$$
and $\mu_0$ denotes the absolutely continuous invariant probability measure.
A result correponding to that of Keller and Nowicki above was obtained 
by Ara\'ujo and Pacifico \cite{AM}  
in more general setting of nonuniformly hyperbolic dynamical systems. 
Melbourne and Nicol \cite{MN} gave an induced scheme approach 
for estimates on the rate functions
in dynamical systems
modelled by Young towers \cite{Y3,Y4} with summable decay of correlations.
All of the results above for nonuniformly hyperbolic dynamical systems 
are obtained under the assumption
of the existence of absolutely continuous invariant probability measures.
But the case that the absolutely continuous
invariant probability measures do not exist has not been considered. 
Also, it is unknown yet neither the criteria to satisfy 
the large deviation principle 
nor the expressions of the rate functions
for nonuniformly hyperbolic dynamical systems. 

The purpose of this paper is to consider the large deviation 
principle for dynamical systems 
from the view point of recurrence times.
We offer a little different description of a tower 
from those already known
to consider a kind of 
specification property
for large deviations estimates.
The topology on a tower in this paper is slightly coarser 
than but almost same as that Young \cite{Y3,Y4} introduced,
and in which we give a sufficient condition on the shape of a 
tower to have a property
that any orbit not recurrent to the base for arbitrarily long time
can be approximated 
by another one recurrent quickly on a tower.
Then a criterion is obtained to ensure the large deviation 
rate functions for abstract dynamical systems.
We show that if a tower satisfies the nonsteep condition
mentioned in the next section,
then the rate function is explicitly represented
by a quantity concerning 
the difference between the metric entropy 
and the Jacobian function. 
The notion of nonsteepness
is independent of the decay rate of the tail.
In fact, it is possible to have the large deviation rate function
for an abstract dynamical system with no absolutely continuous 
invariant probability measures.
On the other hand, there is a tower on which 
we cannot determine the rate function for a dynamical system
although the decay of the tail is exponentially fast.
Some of those examples are provided in the third section.
Combining the argument on large deviations for abstract 
dynamical systems
with a theory for hyperbolic measures \cite{C0,K,KM} 
we establish the
large deviation principle for some class of 
smooth interval maps with nonuniform hyperbolicity.
It is shown the rate function 
coincides with the upper regularization of the free energy
function.
The class of maps for which we can apply the estimates in this paper
contains both of Manneville-Pomeau maps \cite{PSY,PM} and  
Collet-Eckmann unimodal maps \cite{BC,BK,CE,KN,NS,Y2}.
The author thinks that our result is applicable
to a large class of smooth dynamical systems modelled by towers
such as considered in \cite{BLS}.
He also thinks that a theory of multifractal analysis is developed 
from the large deviations estimates of this paper.   
It will be treated in the forthcoming paper \cite{C4}.

\vspace{3mm}

\noindent\textbf{Acknowledgements.}
The author would like to thank to 
Professor Y. Takahashi
for the suggestion on large deviations problem.
He also thanks to 
Professors T. Morita, M. Tsujii, M. Denker and H. Bruin
for valuable comments and encouragement.

\section{Results}

Let $(X, \mathcal{B} , m)$ be 
a finite measure space
and consider a decreasing sequence
$\{ X_k \}_{k=0}^{\infty}$ 
of subsets of $X_0 = X$ with positive measures. 
Assume that 
for each integer $k\ge 0$ 
there is a finite measurable partition 
$\mathcal{I}_k$ of $X_k$ satisfying the following properties:
\begin{enumerate}
\item each $J\in \mathcal{I}_k$ has positive measure;
\item
if $J\in \mathcal{I}_k$ intersects with $X_{k+1}$
then $J\subset X_{k+1}$;
\item
for any $K\in \mathcal{I}_{k+1}$ there is 
$J\in \mathcal{I}_k$ such that $K\subset J$.
\end{enumerate}
Then we call the pair $(Z, \mathcal{A} )$ 
which consists of
the space $$Z := \sqcup_{k=0}^{\infty} X_k \times \{ k \}
\subset X\times \mbox{\boldmath $Z$}^+$$
and
its countable partition 
$$\mathcal{A} :=  
\{ J \times \{ k \} :  J\in \mathcal{I}_k , k=0,1,2,\ldots \}$$
a {\it tower}.
A natural $\sigma$-finite measure $m_Z$ 
on the tower $Z$ is defined by
$m_Z ( A ) 
=\sum_{k=0}^{\infty}m(A_k)$ for 
$A =\sqcup_{k=0}^{\infty} A_k \times \{ k \} \subset Z$
with $A_k \in\mathcal{B} \enspace (k=0,1,2,\ldots ).$  
It is obvious that the measure $m_Z$ is finite iff
$\sum_{k=0}^{\infty} m(X_k)<\infty .$
For a tower $(Z, \mathcal{A} )$  put
$\mathcal{D} := \sqcup_{k=1}^{\infty} \mathcal{D}_k$ where
$$\mathcal{D}_k 
:= \{ J\in \mathcal{I}_{k-1} : J\cap X_{k} =\emptyset \}
\enspace (k=1,2,\ldots) .$$ 
We assume that a bi-nonsingular bijection
$g_J :J \to X$ 
can be taken for each $J\in \mathcal{D} .$
Then we call $T:Z\to Z,$ 
$$T(x,k)  := 
\begin{cases}
(x, k+1)
 & \text{ for } x \in X_{k+1}  , \\
(g_J(x), 0) & \text{ for }
x\in J\in \mathcal{D}_{k+1} ,
\end{cases}$$ 
a {\it tower map}
on $(Z, \mathcal{A} )$ induced by
$\mathcal{G} :=
\{ g_J :J \to X : J\in \mathcal{D} \} .$
Remark that $J\times\{ 0\}\subset Z$ is injectively mapped onto 
$X_0\times\{ 0\}$ by $T^k$ if $J\in\mathcal{D}_k .$
We denote by $k(J)$ the integer $k\ge 1$ such that
$J\in\mathcal{D}_k$, i.e., $T^{k(J)} (J\times\{ 0\}) =X_0\times\{ 0\}.$

In this paper we assume that a tower map $T:Z\to Z$
satisfies
both of the admissibility and 
the bounded distortion conditions  
as below.

\vspace{0.3cm}
\noindent
{\bf The admissibility condition.}
For any sequence 
$\{ J_n \}_{n=0}^{\infty}\subset \mathcal{D}$
there is a unique point $z\in J_0\times \{ 0 \}$ 
such that $T^{k(J_0)+\cdots +k(J_{n-1})} z\in J_{n}\times \{ 0\}$
holds for all $n\ge 1$. 
\vspace{0.3cm}

\noindent 
It follows from the admissibility condition that
for any integer $l\ge 1$
the restriction of $T^l$ to the set
$\displaystyle \cap_{n=0}^{\infty} T^{-nl} 
\left(\sqcup_{A\in\mathcal{K}_l} A \right) $ is 
isomorphic to the full shift of $\sharp\mathcal{K}_l$-symbols
if $\mathcal{K}_l$ is nonempty,
where 
\begin{align*}
\mathcal{K}_l :
= \{ A \in &\vee_{i=0}^{l-1}T^{-i}\mathcal{A} :
A=\cap_{i=0}^{n-1} T^{-{l_i}} (J_i \times \{ 0\}) , \\
&1\le n\le l, \enspace J_0, J_1, \ldots , J_{n-1} \in \mathcal{D}, \\ 
&k(J_0)+ k(J_1) +\cdots +k(J_{n-1}) =l, \enspace l_0=0,\\
& l_i =k(J_0)+\cdots + k(J_{i-1}) \enspace (i=1,\ldots , n-1)
\} ,
\end{align*}
and $\sharp B$ denotes the number of elements of a set $B$.

\vspace{0.3cm}
\noindent
{\bf The bounded distortion condition.}
There are a version ${\rm Jac} (T) >0$ 
of the Radon-Nikodym derivative
$\frac{dm_Z\circ T}{dm_Z} ,$
a constant $D_T \ge 1$ and
a sequence $\{ \varepsilon_k \}_{k=0}^{\infty}$ of positive numbers 
with $\lim_{k\to\infty}\varepsilon_k =0$
such that for any integer $n\ge 1$ and 
$A\in \vee_{i=0}^{n-1} T^{-i}\mathcal{A} ,$
$$
\frac{ {\rm Jac} (T)(z)}{{\rm Jac} (T)(w)}\le e^{\varepsilon_{k(A)}}
\quad\text{ and }\quad 
\frac{\prod_{i=0}^{n-1} {\rm Jac} (T)(T^i(z))}
{\prod_{i=0}^{n-1} {\rm Jac} (T)(T^i(w))} \le D_T$$
hold whenever $z, w\in A ,$ 
where $k(A)=\sharp\{ 1\le j\le n :T^j(A)\subset X_0\times\{ 0 \} \} .$ 

\vspace{0.3cm}

\noindent
We remark that ${\rm Jac} (T) (x,k) =1$ holds 
if $x\not\in J$ for all $J\in\mathcal{D}_k .$
From the bounded distortion condition it follows that
$${D_T}^{-1} 
\le \frac{m_Z(A) \prod_{i=0}^{n-1} {\rm Jac} (T)(T^i(z))}{m_Z (T^n A)} 
\le
D_T$$
holds for any $A\in \vee_{i=0}^{n-1} T^{-i}\mathcal{A}$
and $z\in A$.
\vspace{0.3cm}

We take
$$\mathcal{F} :=\left\{ \psi :Z\to \mbox{\boldmath $R$} :
{\rm a \enspace bounded \enspace function \enspace such \enspace that } 
\lim_{n\to\infty} \text{var}_n (\psi ) = 0 
 \right\} $$
as a class of observable functions on $Z$,
where 
\begin{align*}
\text{var}_n (\psi ) := 
\sup \{ 
\vert \psi (z) - \psi (w) \vert :
 z, w\in A \text{ for some } A\in \vee_{i=0}^{n-1} T^{-i}\mathcal{A} \} 
\end{align*} 
for a function $\psi :Z\to \mbox{\boldmath $R$}$ and $n\ge 1.$
Notice that for any $\psi\in\mathcal{F}$ we have
$\lim_{n\to \infty} {\rm var}_1 (S_n \psi ) / n = 0 ,$
where $$S_n \psi (z) := \sum_{i=0}^{n-1} \psi (T^i (z))$$
for each integer $n\ge 1$. 
In fact, for any $\epsilon >0$ taking $N\ge 1$ 
so that $\text{var}_N (\psi ) \le \varepsilon /2$ we have
\begin{align*}
\vert S_n \psi (z) -  S_n \psi (w)\vert 
&\le \left(\sum_{i=0}^{n-N-1} +\sum_{i=n-N}^{n-1}\right) 
\vert \psi (T^i (z)) - \psi (T^i (w))\vert \\
&\le (n-N-1)\text{var}_N (\psi ) + N \text{var}_1 (\psi ) \\
&\le n\epsilon /2 + n\varepsilon /2 \quad = n\varepsilon 
\end{align*}
whenever $z, w\in A$ for some 
$A\in \vee_{i=0}^{n-1} T^{-i}\mathcal{A}$
and $n\ge 1$ is large enough.

To give large deviations estimates for a tower map
we define the following notion on the shape of a tower.
We say that a tower $(Z, \mathcal{A} )$ is 
{\it nonsteep}, or it satisfies 
{\it the nonsteep condition},
if 
there are sequences
$\{ l_k \}_{k=0}^{\infty}\subset \mbox{\boldmath $N$}$ with 
$\lim_{k\to\infty}  l_k /k  =0$
and
$\{ \gamma_k \}_{k=0}^{\infty} \subset (0, 1)$ with 
$\lim_{k\to\infty} (\log \gamma_k) /k  =0$
such that 
$$\displaystyle  m (J\setminus X_{k+l_k} )  \ge \gamma_k m (J)$$
holds for all $k\ge 0$ and $J\in\mathcal{I}_k .$  
We always assume that the sequences above are monotone,
i.e., $\{ l_k \}_{k=0}^{\infty}$ is nondecreasing
and $\{ \gamma_k\}_{k=0}^{\infty}$
nonincreasing respectively without loss of generality. 
Moreover,
we say that $(Z, \mathcal{A} )$ has 
{\it bounded slope}
if the sequences $\{ l_k \}_{k=0}^{\infty}$ and 
$\{ \gamma_k\}_{k=0}^{\infty}$ above can be taken as constants
respectively.
It is obvious that if the tower has bounded slope,
then the tail decays exponentially fast, 
i.e.,
$ \limsup_{k\to\infty} (\log m(X_k)) / k<0 ,$
and then on which a tower map 
has the exponential decay of correlation
for a function $\psi\in\mathcal{F}$ such that
$\text{var}_n (\psi )$ converges to zero sufficiently fast \cite{Y3}.  

The main result of this paper is the following:
\begin{thm}
Let $(Z, \mathcal{A} )$ be a nonsteep tower and
$T:Z\to Z$ a tower map 
satisfying both of the admissibility and the bounded distortion
conditions.
Then for any 
$\psi\in\mathcal{F}$ 
there exists an upper semicontinuous concave function 
$q_{\psi} :\mbox{\boldmath$R$} \to [-\infty , 0]$ satisfying 
$$ \quad \liminf_{n\to\infty} 
\frac{1}{n} \log m \Bigl( \Big\{ x\in X :
\frac{1}{n} S_n \psi (x , 0) > a \Big\}\Bigr) 
\ge \sup_{t> a} q_{\psi} (t)$$
and 
$$ \quad \limsup_{n\to\infty} 
\frac{1}{n} \log m \Bigl( \Big\{ x\in X :
\frac{1}{n} S_n \psi (x , 0) \ge a \Big\}\Bigr) 
\le \max_{t\ge a} q_{\psi} (t)$$
for all $a\in\mbox{\boldmath $R$} .$
Moreover, the function $q_{\psi}$ above 
can be represented by
\begin{align*}
q_{\psi} (t)
=\lim_{\varepsilon\to 0+} \sup
\Bigl\{ h_{\nu} (T) - \int\log  {\rm Jac} (T) d\nu 
: \nu &\in \mathcal{M}_T  \\
\text{ with } 
\nu ( \sqcup_{k=0}^{K-1} X_k\times\{ k \} )=1 
 &\text{ for some } K\ge 1 \\
\text{ such that  } 
 &| \int \psi d\nu - t | <\varepsilon  \Bigr\} 
\end{align*}
where $\mathcal{M}_T$ denotes the set of 
all $T$-invariant probability measures
on $Z$ and
$h_{\nu} (T)$  
the metric entropy of  $\nu\in \mathcal{M}_T$ for $T$. 
\end{thm}

The theorem above is applicable to large deviations problems for
nonuniformly hyperbolic dynamical systems.
In fact, we obtain a criterion to satisfy the large deviation principle
for smooth interval maps modelled by tower dynamical systems.

Let $I$ be a compact interval of the real line and
$m$ denotes Lebesgue measure on $I$ as a reference measure.
We say that a map $f:I\to I$ is {\it topologically mixing} if
for any nontrivial interval $L\subset I$ 
there is an integer $K\ge 1$
such that $f^KL=I .$  
Let $f:I\to I$ be a $C^2$ map of topologically mixing
and assume that there are a closed subinterval $J$ of $I$, 
a return time function 
$R : J\to\mbox{\boldmath $N$}\cup \{ \infty \}$, i.e.
$f^{R(x)}(x)\in J$ whenever $R(x)<\infty ,$
constants $\lambda >1 ,$   $D \ge 1,$
sequences 
$\{\varepsilon_k\}_{k=0}^{\infty}$ of positive numbers with  
$\lim_{k\to \infty} \varepsilon_k  =0 ,$
$\{ l_k\}_{k=0}^{\infty} \subset\mbox{\boldmath $N$}$
with $\lim_{k\to\infty} l_k /k = 0$
and
 $\{{\gamma}_k\}_{k=0}^{\infty} \subset (0,1)$
with $\lim_{k\to \infty} (\log {\gamma}_k) /k =0$ 
satisfying the following properties:
\begin{enumerate}
\item
if $k\ge 1$ and $V$ is a connected component of $\{ x\in J: R(x) =k \}$, then 
 $$f^kV=J \quad \text{and} \quad 
|(f^k)'(x)|\ge \lambda \enspace (x\in V);$$ 
\item
if $n=k_0+\cdots + k_l\ge 1$ and 
$U_n$ is a connected component of
\begin{align*}
\{ x \in J :R(x)= k_0 , R(f^{k_0}(x))=k_1, \ldots ,
R&(f^{k_0+\cdots +k_{l-2}}(x))= k_{l-1}  \\
\text{and}\enspace &R(f^{k_0+\cdots +k_{l-1}}(x))\ge k_l  
\} ,
\end{align*} 
then 
$m(f^j(U_n)) \le \varepsilon_{n-j}$ for all $0\le j\le n-1.$
\item
if $n=k_0+\cdots + k_l\ge 1$ and 
$V_n$ is a connected component of
$$\{ x \in J :R(x)= k_0 , R(f^{k_0}(x))=k_1, \ldots ,
R(f^{k_0+\cdots +k_{l-1}}(x))=k_l  
\} ,$$ then 
\begin{align*}
\frac{|(f^{k_0})'(y)|}{|(f^{k_0})'(z)|} \le e^{\varepsilon_l}
\quad \text{ and }\quad
\frac{|(f^{n})'(y)|}{|(f^{n})'(z)|} \le D
\end{align*}
hold whenever $y, z\in V_n$;
\item
if $k\ge 1$ and $U$ is a connected component of $\{ x\in J: R(x) >k \}$, 
then
$$
 m \left( \left\{ x\in U : R(x)\le k+l_k \right\} \right) \ge
\gamma_k m(U) .$$
\end{enumerate}
We remark that the function $R$ does not 
necessarily correspond to a first return time on $J$.

We say that $\mu\in \mathcal{M}_f$ is {\it hyperbolic} if
the Lyapunov exponent 
$\lambda (x) :=\limsup_{n\to\infty} (\log |(f^n)'(x) |)/ n $
is positive for $\mu$-almost every $x\in I .$
It follows from the assumptions for the map $f$ 
that hyperbolic measures are dense in $\mathcal{M}_f$.
If $\mu\in \mathcal{M}_f$ is ergodic 
then the Lyapunov exponents coincide with
the constant 
$\displaystyle \lambda_{\mu} (f):= \int \log |f'| d\mu$ 
$\mu$-almost everywhere.
A theory for hyperbolic measures \cite{C0, K, KM} 
asserts the following:
\begin{prop}
Let $\mu\in\mathcal{M}_f$ be ergodic and hyperbolic.
Then, for any continuous function $\varphi : I\to \mbox{\boldmath $R$}$
and $\varepsilon >0$ 
there are integers $k, l\ge 1$
with $(\log l)/k\ge h_{\mu} (f) -\varepsilon$ and
pairwise disjoint compact intervals 
$L_1,L_2, \ldots , L_l$ with $L\subset I$
such that
$L_i\subset L, f^k(L_i) =L$ and 
$L_i$ is injectively mapped to $L$  by $f^k$ for each
$i=1,2,\ldots ,l. $
Moreover, 
$$
| \frac{1}{k} \log\vert (f^k)' (x) \vert 
-\lambda_{\mu} (f) | \le\varepsilon 
\quad \text{ and } \quad 
|\frac{1}{k}S_k \varphi (x) - \int \varphi d\mu |\le
\varepsilon $$
hold whenever $x\in\sqcup_{i=1}^l L_i .$ 
\end{prop}

Now we define {\it the free energy function}
$F: \mathcal{M} \to 
\mbox{\boldmath $R$} \cup \{ -\infty \}$ by
$$F(\mu ) :=
\begin{cases}
 \displaystyle h_{\mu} (f) - \int\log \vert f' \vert d\mu 
 & \text{ for } \mu\in\mathcal{M}_f \text{ hyperbolic,} 
 \\
\displaystyle  -\infty & \text{ otherwise.}
\end{cases}$$
Then combining Proposition 2 with Theorem 1 
we obtain the following:
\begin{thm}[The large deviation principle]
Let $f:I\to I$ be a map satisfying the assumptions above.
Then $f$ satisfies the large deviation principle,
and the rate function $q$ 
coincides with the upper regularization of $F,$ i.e., 
$$q(\mu ) =\inf\{ Q(\mathcal{G} ) :
\mathcal{G} \text{ is a neighborhood of }
\mu  \text{ in } \mathcal{M} \}$$
where $$Q(\mathcal{G} ) :=\sup\{ F(\nu ) : \nu\in\mathcal{G} \} .$$
\end{thm}

\begin{cor}[The Ruelle inequality \cite{Ru}]
For any $\mu\in\mathcal{M}_f$, 
$F(\mu )\le 0$ holds.
\end{cor}

It should be noticed that we need the upper regularization 
for $F$ to get the rate function 
without assuming uniform hyperbolicity of the map $f$,
because the free energy function itself 
may not be upper semiconitinuous 
for a smooth interval map  modelled by a tower dynamical system,
see \cite{BK}.

If $f:I\to I$ is a nonrenormalizable Collet-Eckmann unimodal map $f:I\to I$, 
then a subinterval $J$ can be taken 
with a return time function $R$ satisfying the assumptions above 
so that
the sequences $\{ l_k\}_{k=0}^{\infty}$ and 
$\{ \gamma_k\}_{k=0}^{\infty}$ are constants respectively \cite{Y3},
see also \cite{Ba}.
Then the tower induced from the suspension by the return time function
has bounded slope, and then the tail decays exponentially fast.
It is known that the map $f$ has 
an absolutely continuous invariant probability measure $\mu_0$ and 
the correlation decays exponentially fast for any
continuous function of bounded variation \cite{KN,NS,Y2}.  
It is also known that all of the invariant Borel probability measures 
are hyperbolic for Collet-Eckmann unimodal maps \cite{BK,NS}.

Let $f:I\to I$ be as in Theorem 3 and
$\varphi : I \to \mbox{\boldmath $R$}$
a continuous function.  
Here $\varphi$ is not assumed 
to be of bounded variation. 
Put 
$$ c_{\varphi} := \inf_{x\in I}\liminf_{n\to\infty}
\frac{1}{n} S_n\varphi (x) = 
\min_{\mu\in \mathcal{M}_f } \int\varphi \, d\mu $$ 
and
$$d_{\varphi} := \sup_{x\in I}\limsup_{n\to\infty}
\frac{1}{n} S_n\varphi (x) = 
\max_{\mu\in \mathcal{M}_f } \int\varphi \, d\mu  ,$$
respectively. 
Then the function  
$F_{\varphi} : \mbox{\boldmath$R$}\to [-\infty ,0]$
defined by
$$F_{\varphi} (t) := 
\sup \left\{ F(\mu ) : \int\varphi d\mu = t \right\} $$
is bounded and concave on the interval 
$[c_{\varphi} , d_{\varphi} ] .$
Thus it follows immediately from the theorem that: 
\begin{cor}[The contraction principle]
\begin{align*}
\lim_{n\to\infty}
\frac{1}{n}\log m
\Big(\Big\{ x\in I :
a\le \frac{1}{n}  S_n\varphi (x) \le b \Big\}\Big) 
= \max_{a\le t \le b }  F_{\varphi}(t) 
\end{align*}
holds for any $a , b\in \mbox{\boldmath $R$}$ 
whenever $a\not= d_{\varphi}$ and $b\not= c_{\varphi} .$
\end{cor}
As a consequence we obtain
$$\alpha (\varepsilon ) =
\sup\left\{ F(\mu ) : \vert \int \varphi \, d\mu  
- \int\varphi\, d\mu_0 \vert 
\ge\varepsilon  \right\} $$
for $\alpha$ in the large deviations theorem (\ref{ldp:KN})
for Collet-Eckmann unimodal maps.
The above formula includes the large deviations theorem
because $F(\mu)=0$ holds if and only if $\mu$ is an absolutely continuous
invariant probability measure, i.e., $\mu = \mu_0$ \cite{L1}. 

Another consequence of Theorem 3 follows  
from a general theory on 
large deviations in dynamical systems \cite{T1,T2}.
It is the following:
\begin{cor}[The variational principle of Gibbs type]
The limit 
$$P(\varphi ) :
= \displaystyle \lim_{n\to\infty} 
\frac{1}{n}\log \int \exp S_n\varphi \, dm$$
exists for any continuous function $\varphi : I\to \mbox{\boldmath $R$} .$ 
Moreover, the function 
$P : C(I)\to \mbox{\boldmath $R$}$, the pressure with respect to $m$, 
coincides with the Legendre transform of $-q$,
i.e., 
$$P(\varphi ) 
= \max_{\mu\in \mathcal{M}_f} 
\left\{ q(\mu ) + \int \varphi \, d\mu \right\} 
\quad \text{for all} \quad \varphi\in C(I) , $$
and
$$q(\mu )
= \min_{\varphi\in C(I)}
\left\{ P(\varphi ) - \int \varphi \, d\mu \right\} 
\quad \text{for all} \quad \mu\in \mathcal{M}_f , $$
where $C(I)$ denotes the space of 
the continuous functions on $I$.

\end{cor}

\section{Examples of towers}

In this section we give some of examples of towers.
Throughout this section let 
$X:=(0,1]$ and the measure $m$ on $X$ is Lebesgue measure.
For a sequence $\{a_k\}_{k=0}^{\infty}$ with
 $a_0=1\ge a_1\ge \cdots \ge a_k \ge \cdots >0$ 
setting $X_k := (0, a_k ]$ and 
$\mathcal{I}_k :=\{ X_{k+1} , X_k \setminus X_{k+1} \}$
we obtain a tower $(Z, \mathcal{A} ) $
as in the previous section.

For integers $k$ with $a_{k+1} <a_k$ we define
a linear bijection $g_k =g_{X_k \setminus X_{k+1}} : 
X_k \setminus X_{k+1} \to X $ by
$$g_k (x) =\frac{x-a_{k+1}}{a_k -a_{k+1}} .$$
Then a tower map $T: Z\to Z$ is also defined by 
$$T(x, k) :=
\begin{cases}
 \displaystyle (x , k+1 )  
 & \text{ if } x\in X_{k+1} , 
 \\
\displaystyle  (g_k(x) , 0)  & \text{ if } 
x\in X_{k}\setminus X_{k+1}.
\end{cases}$$
It has no distortion, that is, $D_T=1$ holds in the
bounded distortion condition.

\begin{rmk}
The tower map $T$ defined as above gives a model for 
the countable piecewise linear map $f:[0,1] \to [0,1]$ 
with intermittency introduced originally by Takahashi \cite{T0}:
$$f(x) :=
\begin{cases}
 \displaystyle (x-\beta_1)/(\beta_0-\beta_1)  
 & \text{ for } x\in (\beta_1 , \beta_0 ] , 
 \\
\displaystyle 
\lambda_k (x-\beta_{k+1} ) +\beta_k
 & \text{ for } x\in (\beta_{k+1} , \beta_{k} ] 
 \text{ with } k\ge 1 , \\ 
0 & \text{ for } x=0,
\end{cases}$$
where $\{\beta_k \}_{k=0}^{\infty}$
is a decreasing sequence of positive numbers
with $\beta_0 =1$ and 
$\lambda_k := (\beta_{k-1} -\beta_k )/(\beta_k - \beta_{k+1} )$
for each integer $k\ge 1.$ 
\end{rmk}
\begin{rmk}
By taking another family of functions 
$\mathcal{G} =\{ g_k : a_{k+1} <a_k \}$
it also gives a model on a tower obtained from 
a sequence $\{ a_k\}_{k=0}^{\infty}$ as above
for
a Manneville-Pomeau map, $f(x)=x+x^{1+s} ({\rm mod } 1)$ 
where $0<s<1,$ on the interval $[0,1].$ 
The sequence $\{ a_k\}_{k=0}^{\infty}$ 
corresponds to the preimages of the discontinuity point
of the map.
An estimate on the upper bound is known
for large deviations of this map  \cite{PSY}.
Both of the lower and the upper bounds for large deviations
are obtained from the result of this paper.
\end{rmk}

The stochastic properties of the tower map $T: Z\to Z$ 
are completely determined by the sequence $\{ a_k \}_{k=0}^{\infty} .$
It is well-known that  
the absolutely continuous invariant probability measure
exists for the map $T$ if and only if the sequence is summable, i.e.,
$\sum_{k=0}^{\infty} a_k <\infty .$  
It is also known that the central limit theorem holds 
if $\sum_{k=n}^{\infty} a_k \le Cn^{-\alpha } \enspace (n\ge 1)$ for some 
constants $C\ge 1$ and $\alpha >1 .$
Moreover, if the sequence decays to zero exponentially fast,
i.e., 
$\limsup_{n\to\infty} (\log a_n )/n <0,$  then
so does the correlation function \cite{Y3}.  
However, the large deviations estimates as in  Theorem 1 do not follow
from any conditions mentioned above.
For example,
the sequence $\{ a_k \}_{k=0}^{\infty}$ given by 
$a_k =\exp\{-8^{l+1}\}$
$(8^l\le k< 8^{l+1} , l\ge 0)$
decays to zero exponentially fast,
but the map $T$ does not have a rate function $q_{\psi}$ 
as in Theorem 1
for a locally constant function given by
$\psi (x, k) := 1\enspace (8^l\le k < 2\cdot 8^{l} , l\ge 0 );
=0\enspace (\text{otherwise} ) .$
In fact, since
$$\left\{ x\in X : \frac{1}{8^l} S_{8^l} \psi (x, 0) >
\frac{7}{16} \right\}   =\emptyset$$
for any integer $l\ge 1$, we have
$$\liminf_{n\to\infty} \frac{1}{n}\log
m \left(\left\{ x\in X : \frac{1}{n} S_n \psi (x, 0) >
\frac{7}{16} \right\} \right) = -\infty .$$
On the other hand, since
\begin{align*}
\left\{ x\in X : \frac{1}{2\cdot 8^l}
S_{2\cdot 8^l} \psi (x, 0) \ge
\frac{1}{2} \right\} \supset 
X_{2\cdot 8^l} =
X_{8^l} 
\end{align*}
we have 
$$\limsup_{n\to\infty} \frac{1}{n}\log
m \left(\left\{ x\in X : S_n \psi (x, 0) \ge
\frac{1}{2} \right\} \right) \ge -4.$$
Thus, we cannot take a function $q_{\psi} $ 
to satisfy both of the lower and the upper estimates 
as in Theorem 1.
The tower $(Z , \mathcal{A} )$ obtained from the sequence 
$\{ a_k\}_{k=0}^{\infty}$ as above 
is nonsteep 
if and only if 
$$\lim_{k\to\infty}\frac{1}{k}\log\frac{a_k-a_{k+l_k}}{a_k} =0$$
holds for some sequence $\{ l_k \}_{k=0}^{\infty}$ of positive 
integers
such that $\lim_{k\to\infty} l_k/k=0 .$
Then the rate function $q_{\varphi}$
is given as in Theorem 1 for any 
$\psi\in\mathcal{F} .$  
It is given a typical example of the sequence 
$\{ a_k\}_{k=0}^{\infty}$ for which 
the tower is nonsteep 
by $a_k=(1-p)^{-k}\enspace (k\ge 0)$ with $0<p<1$. 
Then for a function defined by
$\psi (x, k):=1\enspace (k=0) ; =0 \enspace (k\ge 1)$
the rate function $q_{\psi}$
satisfies
$$q_{\varphi} (t)=
\begin{cases}
\displaystyle H(t, 1-t | p, 1-p )  \quad &\text{for} \quad 0\le t\le 1 ,\\
-\infty \quad &\text{otherwise} ,
\end{cases}$$
where $H(t, 1-t | p, 1-p )$ denotes the relative entropy, i.e.,
$$H(t, 1-t | p, 1-p ) := -t\log t- (1-t)\log(1-t) 
+t\log p + (1-t) \log (1-p) .$$ 
This is a classical result on large deviations obtained by
Khinchin \cite{K}.
The tower obtained from this sequence
is not only nonsteep but also having bounded slope. 
In general, the tower obtained from the sequence 
$\{ a_k \}_{k=0}^{\infty}$ has bounded slope if and only if 
$$\frac{a_{k}-a_{k+l}}{ a_k } \ge c \quad (k=0,1,2,\ldots )$$ 
holds for some $l\in\mbox{\boldmath$N$}$ and $c >0 .$  
Another example satisfying the nonsteep condition
is given by    
$a_k = 1/k \enspace ( k\ge 1).$
Then the sequence $\{ a_k\}_{k=0}^{\infty}$ 
is not summable and hence the map $T:Z\to Z$ has no
absolutely continuous invariant probability measures,
nevertheless the large deviations estimates hold 
by Theorem 1. 
It should be noticed that the nonsteep condition
is not necessary for the large deviations estimates.
In fact, the same estimates as in Theorem 1 still hold  
for the tower obtained from the sequence $\{ a_k \}_{k=0}^{\infty}$ 
given by 
$a_k =\exp\{ -8^{2(l+1)}\} \enspace (8^l\le k< 8^{l+1} , m\ge 0),$
although the tower fails the nonsteep condition.
It can be checked that the large deviations estimates 
as in Theorem 1 
valid without the nonsteep condition if the sequence decays
super exponentially fast, i.e., 
$$\displaystyle \lim_{k\to\infty} \frac{1}{k}\log a_k =-\infty$$ 
holds, in general.

\section{Proof of Theorem 1}

Let $(Z,\mathcal{A} )$ be a nonsteep tower, 
and 
$T:Z\to Z$ a tower map satisfying both of the admissibility
and the bounded distortion conditions.
We fix $\psi\in\mathcal{F}$ and  $a\in\mbox{\boldmath $R$}$.
Then the proof of Theorem 1 is devided into two estimates below: 
\begin{enumerate}
\item (The lower estimate)  
$$ \quad \liminf_{n\to\infty} 
\frac{1}{n} \log m \Bigl( \Big\{ x\in X :
\frac{1}{n} S_n \psi (x , 0) > a \Big\}\Bigr) 
\ge \sup_{t> a} q_{\psi} (t) ;$$ \\
\item (The upper estimate)
$$ \quad \limsup_{n\to\infty} 
\frac{1}{n} \log m \Bigl( \Big\{ x\in X :
\frac{1}{n} S_n \psi (x , 0) \ge a \Big\}\Bigr) 
\le \max_{t\ge a} q_{\psi} (t)$$
\end{enumerate}
where 
\begin{align*}
q_{\psi} (t)
=\lim_{\varepsilon\to 0+} \sup
\Bigl\{ h_{\nu} (T) - \int\log  {\rm Jac} (T) d\nu 
: \nu &\in \mathcal{M}_T  \\
\text{ with } 
\nu ( \sqcup_{k=0}^{K-1} X_k\times\{ k \} )=1 
 &\text{ for some } K\ge 1 \\
\text{ such that } 
 &| \int \psi d\nu - t | <\varepsilon  \Bigr\} .
\end{align*}

\noindent
{\bf The lower estimate.}
It is enough to show that
for any  $\nu\in\mathcal{M}_T$ with
$\nu ( \sqcup_{k=0}^{K-1} X_k\times\{ k \} )=1$ 
and $\varepsilon ,\eta >0$ 
the inequality
\begin{align}
m \Big( \Big\{ x\in X : 
| \frac{1}{n} &S_n \psi (x, 0)
-\int \psi d\nu | \le\varepsilon  
\Big\} \Big) \label{ineq:let} 
\\
&\ge
\exp \left\{ n \Big(h_{\nu} (T) -\int\log {\rm Jac} (T) d\nu -\eta 
\Big)  \right\}  \nonumber
\end{align}
holds for any sufficiently large integer $n\ge 1$.
To prove the inequality above, first we assume that $\nu$ is ergodic.
Let 
$$\mathcal{A}_K :=\sqcup_{k=0}^{K-1} \{ J\times\{ k \} : 
J\in\mathcal{I}_k \} \subset \mathcal{A} .$$
Then it is obvious that $\nu (\sqcup_{
A\in\vee_{i=0}^{n-1} T^{-i}\mathcal{A}_K} A)= 1$
since $\nu\in\mathcal{M}_T .$ 
Also, since $\psi \in \mathcal{F} $
$$| S_n\psi (z)-
S_n\psi (w) |\le n\varepsilon /8$$
holds whenever 
$z, w\in A$ for some $A\in\vee_{i=0}^{n-1} T^{-i}\mathcal{A}_K$
and $n\ge 1$ is large.
Let 
\begin{align*}
\mathcal{B}_n :=
\Big\{ A\in&\vee_{i=0}^{n-1} T^{-i}\mathcal{A}_K : 
\nu (A) \le \exp\{ -n(h_{\nu} (T)-\eta /8 ) \} , \\
&| \sum_{i=0}^{n-1} \log {\rm Jac}(T) (T^i(w_A)) 
-n\int \log {\rm Jac}(T) d\nu |\le n\eta /8 \\
&\text{and }\enspace 
| S_n\psi (z_A) 
-n\int \psi d\nu |\le n\varepsilon / 8
\enspace \text{ for some } z_A, w_A \in A 
\Bigr\}
\end{align*}
for each integer $n \ge 1.$
By the Birkhoff ergodic theorem and 
the Shannon-McMillan-Breimann theorem 
$\nu (\sqcup_{A\in \mathcal{B}_n} A ) \ge 1/2$
holds, 
and hence the number of elements in $\mathcal{B}_n$
is not smaller than $e^{ n(h_{\nu}(T)-\eta /8 ) } /2$
for large $n\ge 1.$
Then we can choose integers $k$ and $l$ with $0\le k, l \le K-1$
so that the set
$$\mathcal{B}_{n, k,l} := 
\{ A\in \mathcal{B}_n : A\subset X_k\times \{ k \} ,
\enspace
T^{n+l}A = X_0 \times \{0 \} \}$$ 
contains at least $e^{ n(h_{\nu}(T)-\eta /4 ) }
\enspace (\le e^{ n(h_{\nu}(T)-\eta /8 ) } /(2K^2) )$
elements.
For large $n\ge 1$ and $A\in\mathcal{B}_n$
\begin{align*}
m_Z (A) 
&\ge {D_T}^{-1} m_Z (T^nA) 
\prod_{i=0}^{n-1} {\rm Jac} (T) (T^i(w_A))^{-1} \\
&\ge {D_T}^{-1} \min_{B\in\mathcal{A}_K} m_Z (B) 
\exp\left\{ -n\left(\int\log{\rm Jac}(T) d\nu + \eta /8\right) \right\} \\
&\ge 
\exp\left\{ -n\left(\int\log{\rm Jac}(T) d\nu + \eta /4\right) \right\} 
\end{align*}
holds.
Then we have 
\begin{align*} \label{ineq:tri1}
\sum_{A\in\mathcal{B}_{n,k,l}} m_Z(A)  
&\ge \sharp \mathcal{B}_{n,k,l}\cdot
\min_{A\in\mathcal{B}_{n,k,l}}m_Z(A) \\
&\ge  
\exp\Big\{ n\Big(h_{\nu} (T) - \int\log {\rm Jac}(T) d\nu 
- \eta /2\Big) \Big\} . \nonumber
\end{align*}
Take $A^* \in \vee_{i=0}^{n+k-1} T^{-i}\mathcal{A}_K $
such that $T^kA^* =A$ for each $A\in \mathcal{B}_{n, k,l} .$
Then it is obvious that $A^* \subset X_0\times \{ 0 \} $
and $T^{n+k+l}A^* = X_0 \times \{ 0\} .$
Let
$$\mathcal{B}_{n, k,l}^* := 
\{ A^* : 
A\in \mathcal{B}_{n, k,l} \}  .$$ 
Then since $m_Z(A^*) =m_Z(A)$ for each $A\in \mathcal{B}_{n, k,l} $ 
we have 
\begin{equation*}\label{ineq:tri2}
\sum_{A^*\in\mathcal{B}^*_{n,k,l}} m_Z(A^*)
=\sum_{A\in\mathcal{B}_{n,k,l}} m_Z(A).
\end{equation*}
Moreover,
for any $A\in\mathcal{B}_{n,k,l} $ and $z\in A$
there is $z^*\in A^*$ such that $T^k(z^*)=z ,$
and then 
\begin{align*} 
| &S_n\psi (z^*) -n\int\psi d\nu | \\ 
&\le |S_n\psi (z^*)-  S_n \psi (z) | \nonumber
+ | S_n \psi (z) - S_n \psi (z_A) | + | S_n\psi (z_A) -n\int\psi d\nu | \\
&\le 2K\sup_{w\in Z} |\psi (w) | + n\varepsilon /8 + n\varepsilon /8   
\qquad \le n\varepsilon /2 . \nonumber
\end{align*}
As a consequence we obtain
\begin{align*}
m\Big( \Big\{ &x\in X : 
| \frac{1}{n}S_n\psi (x,0) - \int\psi d\nu |\le  \varepsilon /2  
 \Big\} \Big) \\
&= m_Z \Big( \Big\{ z^* \in X_0\times\{ 0 \} : 
| S_n\psi (z^*) -
n \int\psi d\nu |\le  n\varepsilon /2  
 \Big\} \Big) \\
&\ge\sum_{A^*\in\mathcal{B}^*_{n,k,l}} m_Z(A^*)
\qquad =\sum_{A\in\mathcal{B}_{n,k,l}} m_Z(A) \\
&\ge \exp\Big\{ n\Big( h_{\nu} (T) - \int\log {\rm Jac} (T) d\nu
 - \eta /2\Big) \Big\} . 
\end{align*}
The inequality (\ref{ineq:let}) is proved for the case 
$\nu$ is ergodic.
For $\nu\in\mathcal{M}_T$ not ergodic 
take a linear combination
$\nu ' =\alpha_1 \nu_1 +\cdots + \alpha_p\nu_p$
of ergodic $T$-invariant probality measures 
$\nu_1 ,\ldots , \nu_p$
supported on $\sqcup_{k=0}^{K-1} X_k\times \{ k \}$
such that
\begin{align*}
| h_{\nu}(T) -  h_{\nu '}(T) |\le \eta /8, \quad
| \int\log{\rm Jac} (T)d\nu -\int\log{\rm Jac} (T)d\nu 
|\le \eta /8 
\end{align*}
and
$$| \int \psi d\nu -\int \psi d\nu ' |\le \varepsilon /4 .$$
For large $n\ge 1$ and $q=1,\ldots , p$ put
$n_q :=[n\alpha_q]$ where $[\cdot ]$ denotes the Gauss' symbol.
Applying the above argument  for $\nu_q$
we can take integers $k_q, l_q$ with $0\le k_q ,l_q \le K-1$
and $\mathcal{B}_n^* (q) 
\subset\vee_{i=0}^{n_q+k_q-1} T^{-i}\mathcal{A}_K$
which consists of at least 
$\exp\{ n\alpha_q (h_{\nu_q}(T) -\eta/4)\}$ elements
$A$ such that:
\begin{enumerate}
\item $A\subset X_0\times \{ 0 \} \enspace\text{ and }\enspace 
T^{n_q+k_q+l_q} A = X_0\times\{ 0\}$;
\item $\displaystyle  m_Z(A)\ge 
\exp\left\{ -n\alpha_q
\left(\int\log{\rm Jac}(T) d\nu_q + \eta /8\right) \right\} ;$
\item 
$\displaystyle | S_{n_q}\psi (z) -n_q \int\psi d\nu_q |  
\le n_q\varepsilon /2$ 
holds whenever $z\in A .$
\end{enumerate}
Then for any $z\in A$ with $A\in \mathcal{B}_n^* (q)$ we have
\begin{align*}
\prod_{i=0}^{r_q(n)-1}{\rm Jac}&(T)(T^i(z))
\ge {D_T}^{-1} \frac{m_Z (X_0\times\{ 0 \})}{m_Z(A)} \\
&\ge  {D_T}^{-1} m (X_0)
\exp\left\{ -n_q\left(\int\log{\rm Jac}(T) d\nu_q 
+ \eta /8\right) \right\} \\
&\ge \exp\left\{ -n\alpha_q\left(\int\log{\rm Jac}(T) d\nu_q 
+ \eta /4\right) \right\} 
\end{align*}
where $r_q(n) = n_q +k_q + l_q .$ 
Let
\begin{align*}
\mathcal{B}_n^* := \{ 
B=\cap_{q=1}^{p} T^{-s_{q-1}(n)}B_q  :
B_q\in \mathcal{B}_n^* (q) \text{ for all } q=1,\ldots , p \} 
\end{align*}
where $s_0(n)=0$ and
$s_q(n) =r_1(n)+\cdots + r_q(n) $
for $q=1,\ldots , p .$
Then the number of elements of $\mathcal{B}_n^*$
is not smaller than  
$\displaystyle e^{ n (h_{\nu '}(T) -\eta/4) }  .$
For each $B\in \mathcal{B}_n^*$ take $B_q\in \mathcal{B}_n^* (q)$
for $q=1,\ldots ,p$ with  
$B= \cap_{q=1}^{p} T^{-s_{q-1}(n)}B_q$ 
and $z\in B .$
Then we have
$T^{s_{q-1}(n)}(z)\in B_{q}$
for each $q=1,\ldots , p$, and hence
\begin{align*}
m_Z (B) 
&\ge D_T^{-1}m_Z(X_0\times \{ 0\} )
\prod_{i=0}^{s_q(n)-1}{\rm Jac}(T)(T^i(z))  \\
&= D_T^{-1}m(X_0)
\prod_{q=1}^p \prod_{i=0}^{r_q(n)-1}{\rm Jac}(T)(T^i(T^{s_{q-1}(n)}(z))) \\
&\ge D_T^{-1}m(X_0)
\exp\left\{ -n \sum_{q=1}^p \alpha_q\left(\int\log{\rm Jac}(T) d\nu_q 
+ \eta /4\right) \right\}  \\
&\ge
 \exp\left\{ -n \left(\int\log{\rm Jac}(T) d\nu ' 
+ \eta /2\right) \right\} .
\end{align*}
Therefore,
\begin{align*}
\sum_{B\in \mathcal{B}_n^*} m_Z (B)  
&\ge \sharp \mathcal{B}_n^* \cdot \min_{B\in \mathcal{B}_n^*} m_Z (B) \\
&\ge \exp\left\{ n \left( h_{\nu '}(T) - \int\log {\rm Jac}(T) d\nu ' 
- 3\eta /4 \right) \right\} \nonumber \\
&\ge \exp\left\{ n \left( h_{\nu }(T) - \int\log {\rm Jac}(T) d\nu  
- \eta /4 \right) \right\} .
\nonumber
\end{align*}
Moreover, for any $z\in B$ with $B\in\mathcal{B}_n^*$,
\begin{align*}
| S_n &\psi (z) -n \int\psi d\nu | 
\le
| S_n \psi (z) -n \int\psi d\nu ' |
+ | n \int\psi d\nu '  - n \int\psi d\nu |  \\  
\le &
\sum_{q=1}^p \Big(| S_{n_q}\psi (z) -n_q \int\psi d\nu_q |  
+ 2(k_q+r_q)\sup_{w\in Z} |\psi (w)|\Big) \nonumber \\
+ &n |  \int\psi d\nu '  - \int\psi d\nu_q |  \nonumber \\
\le & \sum_{q=1}^p (n_q\varepsilon /2 +  4K \sup_{w\in Z} |\psi (w)|)
+  n\varepsilon /4 \quad
\le n \varepsilon \nonumber .
\end{align*} 
From the estimates above 
we conclude
\begin{align*}
m\Big( \Big\{ x\in &X : 
| \frac{1}{n}S_n\psi (x,0) - \int\psi d\nu |\le  \varepsilon   
 \Big\} \Big) \\
&= m_Z \Big( \Big\{ z \in X_0\times\{ 0 \} : 
| S_n\psi (z) -
n \int\psi d\nu |\le  n\varepsilon
 \Big\} \Big) \\
&\ge \sum_{B\in \mathcal{B}_n^*} m_Z (B) \\ 
&\ge \exp\left\{ n (h_{\nu }(T) 
 -n \left(\int\log{\rm Jac}(T) d\nu -\eta \right) \right\} .
\end{align*}
The inequality (\ref{ineq:let}) is obtained for $\nu \in\mathcal{M}_T$ 
even if it is not ergodic.
It finishes the proof of the lower estimate.

\begin{rmk}
The nonsteep condition is not needed for
the lower estimate in the above proof.
\end{rmk}

\noindent
{\bf The upper estimate.}
To obtain the upper estimate 
we need a variational principle 
for dynamical systems of bounded distortion as below.
Let $(Y, \mathcal{B} , m)$ be a finite measure space and
$Y_1 , \ldots , Y_l \in\mathcal{B} $ 
pairwise disjoint subsets of $Y$ with 
positive measures.
We consider a measurable map $g: \sqcup_{j=1}^l Y_j \to Y$ 
satisfying the following properties:
\begin{enumerate}
\item for each $j=1,\ldots , l, $
$g_j :=g|_{Y_j} : Y_j\to Y$ is a bi-nonsingular bijection;
\item for any sequence $\{  a_j \}_{j=0}^{\infty}$ with 
$a_j\in \{ 1,\ldots ,l \} \enspace (j=0,1,2,\ldots )$ 
the set $\cap_{n=0}^{\infty} Y_{a_0\ldots a_{n-1}}$ 
consists of a single point,
where $Y_{a_0\cdots a_{n-1}} :
= Y_{a_0}\cap g^{-1}Y_{a_1}\cap\cdots\cap g^{-(n-1)}Y_{a_{n-1}}$;
\item there are a version ${\rm Jac} (g) >0$ 
of the Radon-Nikodym derivative $\frac{dm\circ g}{dm}$ 
with 
$\displaystyle \lim_{n\to\infty}\sup_{J\in W_n}\sup_{x, y\in Y_J}
| \log {\rm Jac}(g)(x)-\log {\rm Jac}(g)(y) | =0$
and a distortion constant $C\ge 1$ such that
for any integer $n\ge 1$ and
$J \in W_n$
$$\frac{\prod_{i=0}^{n-1} {\rm Jac} (g) (g^i(x)) }
{\prod_{i=0}^{n-1} {\rm Jac} (g) (g^i(y)) }
\le C$$
holds whenever $x, y\in Y_J,$
where 
$W_n$ stands for the set of the words 
$J=(a_0 \cdots a_{n-1})$
of length $n$ with 
$a_i\in \{ 1,2,\ldots , l \}$ for each $i=0,1,\ldots ,n-1.$
\end{enumerate}
Then we call $g:\varLambda \to\varLambda$ a 
{\it finite Markov system} induced by $( Y, \{ Y_i\}_{i=1}^l , g) ,$
where $\varLambda := \cap_{n=0}^{\infty} g^{-n} (\sqcup_{j=1}^l Y_j) .$
It is isomorphic to a full shift of $l$-symbols, 
and the space of the probability measures supported on $\varLambda$
is compact.
The following lemma is obtained from a standard argument 
on the variatinal principle for pressure \cite{W}.

\begin{lem}
A finite Markov system $g:\varLambda \to \varLambda$ 
induced by $(Y, \{ Y_j\}_{j=1}^l , g)$
with a distortion constant $C\ge 1$ has an
invariant probability measure $\mu$ on $\varLambda$ such that
$$h_{\mu} (g) -\int\log {\rm Jac} (g) d\mu \ge
\log \sum_{p=1}^l m(Y_p) - \log  m(Y) -  \log C .$$
\end{lem}

\noindent
{\it Proof of Lemma 7.}
For each $J=(a_0\ldots a_{n-1})\in W_n$  and 
$ p= 1,\ldots , l$
let $Y_{Jp} := Y_J\cap g^{-n}Y_p .$ 
Then since
\begin{align*}
\frac{m(Y_{Jp})}{m(Y_{J})} 
\ge \frac{\prod_{i=0}^{n-1}\inf_{x \in Y_{a_0\ldots a_{n-1}p}}
{\rm Jac} (g)(g^i(x))^{-1} m(Y_p)}
{\prod_{i=0}^{n-1}\sup_{y\in Y_{a_0\ldots a_{n-1}}}
{\rm Jac} (g)(g^i(y))^{-1} m(Y)} 
\ge C^{-1} \frac{m(Y_p)}{m(Y)} 
\end{align*}
for each $p=1,\ldots , l,$
we have 
$$\frac{\sum_{p=1}^{l} m(Y_{Jp})}{m(Y_{J})} \ge 
C^{-1} \frac{\sum_{p=1}^{l} m(Y_p)}{m(Y)}  
,$$
and hence
\begin{align*}
\sum_{J\in W_n}  m(Y_J)  
=&\sum_{a_0, \ldots , a_{n-1} =1}^{l}  m(Y_{a_0\ldots a_{n-1}})  \\
\ge &\left\{ C^{-1} \frac{\sum_{p=1}^l m(Y_p) }{m(Y)} \right\}
\sum_{a_0, \ldots , a_{n-2} =1}^{l}  m(Y_{a_0\ldots a_{n-2}}) \\  
\cdots \ge &\left\{
C^{-1} \frac{ \sum_{p=1}^l m(Y_p) }{m(Y)} \right\}^{n-1}
\sum_{a_0 =1}^{l}  m(Y_{a_0}) \\
= &\left\{ 
C^{-1} \frac{ \sum_{p=1}^l m(Y_p)}{m(Y)} \right\}^{n}
Cm(Y) .
\end{align*}
Thus we obtain
\begin{align}
\liminf_{n\to\infty}
\frac{1}{n} \log\sum_{J\in W_n} m(Y_J)
\ge \log \sum_{p=1}^l m(Y_p)  -\log m(Y) - \log C 
 .\label{eq:fms1}
\end{align}
On the other hand, taking $x_J\in Y_J\cap\varLambda$ 
for each $J\in W_n $ 
we obtain a sequence of probability measures supported on $\varLambda$
by
$\displaystyle \mu_n :=\frac{1}{Z_n} \sum_{J\in W_n} m(Y_J) \delta_{x_J}^n $
for each $n\ge 1$,
where $Z_n := \sum_{J\in W_n} m(Y_J)$ and 
$\delta_{x_J}^n := (\delta_{x_J} +\delta_{g(x_J)} + 
\cdots + \delta_{g^{n-1}(x_J)} )/n .$
Then an accumulation point  
$\mu$ of the sequence $\{ {\mu}_n \}_{n=1}^{\infty}$
is a $g$-invariant probability measure 
supported on $\varLambda .$
Since
\begin{align*}
\log &\sum_{J\in W_n} m(Y_J) 
=\log Z_n \\
=
&\sum_{J\in W_n} \left( \frac{m(Y_J)}{\sum_{J'\in W_n}m(Y_{J'})}\right)
\left\{ 
 -\log \left( \frac{m(Y_J)}{\sum_{J'\in W_n}m(Y_{J'})} \right) 
 +\log m(Y_J)   \right\} \\
=
&\sum_{J\in W_n} \mu_n (Y_J) 
\left\{ - \log \mu_n (Y_J) + \log m(Y_J)   \right\} \\
\le
&\sum_{J\in W_n} \mu_n (Y_J) 
\left\{  - \log \mu_n (Y_J) + 
\log \left( C m(Y)\prod_{i=0}^{n-1}
\text{Jac} (g)(g^i(x_J))^{-1} \right)  \right\} \\
=&-\sum_{J\in W_n} \mu_n (Y_J)\log \mu_n (Y_J) 
- n\int\log {\rm Jac} (g) d\mu_{n} +\log \{ C m(Y) \} ,
\end{align*}
we have
\begin{equation}
\limsup_{n\to\infty}\frac{1}{n} \log\sum_{J\in W_n} m(Y_J)
\le h_{\mu} (g) - \int\log {\rm Jac} (g) d\mu . \label{eq:fms2}
\end{equation}
Combining (\ref{eq:fms1}) and (\ref{eq:fms2})
we obtain
$$h_{\mu} (g) -\int\log {\rm Jac} (g) d\mu \ge
\sum_{j=1}^l \log m(Y_j) - \log  m(Y) -  \log C .$$
The lemma is proved. 

\vspace{0.3cm}
Now we show the upper estimate.
Take the monotone sequences 
$\{ l_k \}_{k=0}^{\infty} \subset \mbox{\boldmath $N$}$
and  $\{ {\gamma}_k \}_{k=0}^{\infty} \subset (0, 1)$
as in the definition of the nonsteep condition for $T$.
Put
\begin{align*} 
-\beta :=\limsup_{n\to\infty} 
\frac{1}{n} \log m\Bigl( \Bigl\{ x\in X :
 \frac{1}{n} S_n \psi  (x , 0)
\ge a   \Bigr\} 
\Bigr)  .
\end{align*}
If $\beta =\infty$ then nothing has to be shown,
and so we assume that $\beta <\infty .$
We show that for any $\varepsilon  , \eta >0$
there is $\nu\in\mathcal{M}_T$
with $\nu ( \sqcup_{k=0}^{K-1} X_k\times\{ k \} ) =1$ 
for some $K\ge 1$ such that
$$\int \psi d\nu > a-\epsilon
\quad\text{ and }\quad 
h_{\nu} (T) - \int\log {\rm Jac} \it (T) d\nu \ge -(\beta +\eta ) $$
hold.
Take an arbitrarily large integer $N\ge 1$ such that 
$$e^{ -N(\beta +\eta /4) } \le
 m\Bigl( \Bigl\{ x\in X :
\frac{1}{N} S_N \psi (x , 0) \ge a  \Bigr\}$$ 
and $$|S_N \psi  (z) - S_N \psi  (w) |\le \epsilon N/8$$
whenever $z, w\in A$  for some 
$A\in \vee_{i=0}^{N-1} T^{-i}\mathcal{A} .$
Let
\begin{align*}
\mathcal{B}_N :=
\Bigl\{
A&\in{\vee}_{i=0}^{N-1}T^{-i}\mathcal{A} : 
A\subset X_0\times \{  0 \} , \\ 
&\frac{1}{N} S_N \psi (z_A) \ge a 
\text{ for some } z_A\in A \Bigr\} .
\end{align*}
Then it is obvious that
$$\sum_{A\in \mathcal{B}_N } m_Z(A)
\ge m\Bigl( \Bigl\{ x\in X :
\frac{1}{N} S_N \psi (x , 0) \ge a  \Bigr\}
\ge e^{-N(\beta + \eta /4)} .$$
For each $A\in\mathcal{B}_N$ take $0\le k_A \le N-1$
and $J_A\in\mathcal{I}_{k_A}$ 
 such that
$$T^{N-1}A = J_A \times \{ k_A \} ,$$
and let 
$$J_{A, j} := (J_A\cap X_{k_A+j-1} )\setminus X_{k_A+j} \quad
(j=1,\ldots l_N ).$$
Then from the monotonicity of the sequences 
$\{l_k\}_{k=0}^{\infty}$ and $\{\gamma_k\}_{k=0}^{\infty}$
in the assumption on the nonsteepness
we have
\begin{align*}
\sum_{j=1}^{l_N} m (J_{A,j} ) &= m(J_A\setminus X_{k_A+l_N} ) 
\ge m(J_A\setminus X_{k_A+ l_{k_A}} ) \\  
&\ge \gamma_{k_A} m(J_A) 
\qquad\quad\ge \gamma_N m(J_A) .
\end{align*}
For  $j=1,2,\ldots , l_N$
we  set
$$\mathcal{B}_{N,j} := 
\left\{ B\in  \vee_{i=0}^{N+j-1} T^{-i}\mathcal{A} :
B\subset A_j , A\in\mathcal{B}_N
\right\} ,$$
where $A_j :=A\cap T^{-N+1} (J_{A, j}\times \{ k_A \} ) $
for each $A\in\mathcal{B}_N .$
Then since 
\begin{align*}
\frac{\sum_{j=1}^{l_N} m_Z (A_j)}{m_Z(A)} &\ge
{D_T}^{-1}\frac{\sum_{j=1}^{l_N} m_Z (J_{A,j} \times \{ k_A\} )}
{m_Z (J_{A} \times \{ k_A\} )  } \\
&= 
{D_T}^{-1}\frac{\sum_{j=1}^{l_N} m (J_{A,j} )}
{m (J_{A} )  } \quad \ge
{D_T}^{-1} \gamma_N 
\end{align*}
for each $A\in \mathcal{B}_N ,$
we have 
\begin{align*}
\sum_{j=1}^{l_N}\sum_{ B\in \mathcal{B}_{N,j}} m_Z(B) 
&=\sum_{A\in\mathcal{B}_N}\sum_{j=1}^{l_N} m_Z(A_j) \\
&\ge {D_T}^{-1}\gamma_N  \sum_{A\in \mathcal{B}_{N}} m_Z(A) \\
&\ge  {D_T}^{-1}\gamma_N e^{ -N(\beta +\eta /4)} ,
\end{align*}
and hence,
\begin{align*}
\sum_{B\in \mathcal{B}_{N,j_N}} m_Z(B) \ge
{D_T}^{-1}{l_N}^{-1} \gamma_N  e^{ -N(\beta +\eta /4)} 
\end{align*}
holds for some $1\le j_N\le l_N .$
Then $( X_0\times \{  0\} , 
\mathcal{B}_{N, j_N} , T^{N+j_N}|_{X_0\times \{ 0 \}} )$ induces  
a finite Markov system on $X_0\times \{ 0 \}$.
Set $K:=N+j_N.$  Then we obtain
a $T^{K}$-invariant probability measure
$\mu$ on $\varLambda :=\cap_{l=0}^{\infty} T^{-lK} 
(\sqcup_{B\in \mathcal{B}_{N, j_N}} B )\subset X_0\times \{ 0 \}$
such that 
\begin{align*}
h_{\mu} (T^{K}) &-\int\log {\rm Jac} (T^{K} )d\mu \\
&\ge \log\sum_{B\in\mathcal{B}_{N, j_N}} m_Z(B) -\log m_Z(X_0\times \{ 0 \} )
-\log D_T
\end{align*}
by Lemma 7.
Then $\nu := \frac{1}{K}\sum_{i=0}^{K-1} \mu\circ T^{-i}$
is a $T$-invariant probability measure 
satisfying
$\nu (\sqcup_{k=0}^{K-1} ( X_k\times \{ k \} ))=1 .$
Moreover,
\begin{align*}
h_{\nu} (T) &-\int\log {\rm Jac} (T )d\nu 
=\frac{1}{K} \Big\{ h_{\mu} (T^{K}) 
-\int\log {\rm Jac} (T^{K} ) d\mu\Big\} \\
\ge &\frac{1}{K} \Big\{
\log\sum_{B\in\mathcal{B}_{N, j_N}} 
m_Z(B) -
 \log m_Z(X_0\times \{ 0 \} ) - \log D_T \Big\} \\
\ge &\frac{1}{K} \log\sum_{B\in\mathcal{B}_{N, j_N}} m_Z(B) -
\eta /4  \\
\ge &\frac{1}{K} \log \{
{D_T}^{-1}{l_N}^{-1} {\gamma}_N   e^{ -N(\beta +\eta /4)} \}
 - \eta /4  \\
\ge &-\frac{N}{N+j_N} (\beta +\eta /4) 
- \frac{1}{N} \log  {D_T}  - \frac{1}{N} \log {l_N}
+ \frac{1}{N} \log \gamma_N
- \eta /4   \\
\ge &-(\beta +\eta )
\end{align*}
holds if $N$ is large.
Furthermore, for any 
$z\in B$ with $B\in \mathcal{B}_{N, j_N} $ 
we can take $A\in \mathcal{B}_N$ with $B\subset A$
and $z_A\in A$ such that
$S_N\psi (z_A) / N \ge a.$
Then for a large integer $N\ge 1$ we have
\begin{align*}
S_K\psi (z_A) 
&\ge S_N\psi (z_A) -j_N\sup_{w\in Z} |\psi (w) | \\
&\ge Na - l_N \sup_{w\in Z} |\psi (w)| \\
&\ge (N+j_N)(a -\varepsilon /4 ) \quad
= K(a-\varepsilon /4) ,
\end{align*}
and
\begin{align*}
| S_K\psi (z) 
- S_K\psi (z_A) |
&\le | S_N\psi (z) 
- S_N\psi (z_A) | \\
&+ | S_{j_N} \psi (T^N(z)) 
- S_{j_N} \psi (T^N (z_A)) |  \\
&\le N \varepsilon / 8 + 2j_N \sup_{w\in Z} |\psi (w)| \\
&\le N \varepsilon / 8 + 2l_N \sup_{w\in Z} |\psi (w)| \\
&\le N\varepsilon /4 \quad \le K\varepsilon /4.
\end{align*}
Then,
\begin{align*}
S_K \psi (z)
&= S_K\psi (z_A) + 
\left( S_K\psi (z) - S_K\psi (z_A) \right) \\
&\ge K(a-\varepsilon /4)  -  K\varepsilon /4 \quad
=  K(a-\epsilon /2) .
\end{align*}
Since $\nu\in\mathcal{M}_T$ and it is  supported on  
$\cap_{l=0}^{\infty} T^{-lK} 
(\sqcup_{B\in \mathcal{B}_{N, j_N}}
\sqcup_{i=0}^{K-1} T^iB )$
we have 
$$\int \psi d\nu 
\ge a -\epsilon /2 > a -\epsilon .$$
Thus we obtain the upper estimate.
This completes the proof of Theorem 1.

\begin{rmk}
To obtain the large deviations estimates as in Theorem 1, 
we can relax the bounded distortion condition to weaker one
by replacing the constant $D_T\ge 1$ with
a sequence $\{ D_n \}_{n=1}^{\infty}$ of positive numbers 
satisfying $\lim_{n\to \infty} (\log D_n) / n = 0$
such that 
$$\frac{\prod_{i=0}^{n-1} {\rm Jac} (T)(T^i(z))}
{\prod_{i=0}^{n-1} {\rm Jac} (T)(T^i(w))} \le D_n$$
holds whenever 
$z, w\in A$ for some $A\in \vee_{i=0}^{n-1} T^{-i}\mathcal{A}$ and
$n\ge 1 .$
A thermodynamic formalism for dynamical systems satisfying
this weak bounded distortion condition  
has been studied  by Yuri \cite{Yu1}.
\end{rmk}

\section{Proof of Theorem 3}

Throughout this section, $I$ denotes a compact interval of
the real line and $m$ Lebesgue measure.
Let $f: I\to I$ be a topological mixing $C^2$ map satisfying the assumptions
stated in Section 2.
We notice that the weak* topology on the space 
$\mathcal{M}$ of the probability measures
is generated by open sets 
$\mathcal{G}$ of the form
$$\mathcal{G} :=\Big\{ \mu\in\mathcal{M} : 
\max_{i=1,\ldots ,l} \Big| \int \varphi_i d\mu - \alpha_i \Big| < \varepsilon
\Big\}$$
for some $\varphi_1 ,\ldots , \varphi_l \in C(I) ,$  
$\alpha_1, \ldots , \alpha_l\in\mbox{\boldmath$R$}$
and $\varepsilon >0,$
where $C(I)$ denotes the space of the continuous functions on $I$.  
Thus, the proof of Theorem 3  is reduced to the estimates
on the level 1 large deviations for given
$\varphi \in C(I)$ and $a\in \mbox{\boldmath$R$}$
as follows:
\begin{enumerate}
\item (The lower estimate)  
\begin{align*} \liminf_{n\to\infty} 
\frac{1}{n} \log m \Bigl( \Big\{ x\in I :
\frac{1}{n} S_n \varphi (x) > a \Big\}\Bigr)  
\ge \sup \Big\{ q (\mu ) : \int \varphi d\mu > a \Big\}   ;
\end{align*} \\
\item (The upper estimate)
\begin{align*}
\limsup_{n\to\infty} 
\frac{1}{n} \log m \Bigl( \Big\{ x\in I :
\frac{1}{n} S_n \varphi (x) \ge a \Big\}\Bigr) 
\le \max \Big\{ q (\mu ) : \int \varphi d\mu \ge a \Big\} ,
\end{align*}
\end{enumerate}
where 
$q$ is the upper regularization of $F$,
$$F(\mu ) :=
\begin{cases}
 \displaystyle h_{\mu} (f) - \int\log \vert f' \vert d\mu 
 & \text{ for } \mu\in\mathcal{M}_f \text{ hyperbolic,} 
 \\
\displaystyle  -\infty & \text{ otherwise.}
\end{cases}$$

\noindent
{\bf The lower estimate.}
To obtain the lower estimate of the rate function
we show that 
\begin{align}
\liminf_{n\to\infty}
\frac{1}{n}\log m \left( \left\{ x\in I :
\frac{1}{n} S_n\varphi (x) > a \right\} \right)  
\ge h_{\mu} (f) - \int \log |f' | d\mu \label{ineq:les} 
\end{align}
holds for any
$\mu\in\mathcal{M}_f$ hyperbolic with $\int \varphi d\mu > a .$
First we assume that $\mu$ is ergodic.
Taking $\varepsilon >0$ small enough so that 
$\int \varphi d\mu > a + \varepsilon .$
Then there are integers $k, l\ge 1$
with $(\log l)/k\ge h_{\mu}(f) -\varepsilon$ and
pairwise disjoint compact intervals 
$L_1,L_2, \ldots , L_l$ with $L\subset I$
such that:
$L_i\subset L , f^k(L_i)=L$ and 
$f^k\mid_{L_i} :L_i \to L$ is injective on $L_i$
$(i=1,2,\ldots l)$;
$$| \frac{1}{k}\log\vert (f^k)' (x) \vert 
- \int \log |f' | d\mu | \le\varepsilon
\quad\text{ and }\quad
\frac{1}{k} S_k\varphi (x) \ge \int \varphi d\mu -\epsilon
> a  $$
whenever $x\in\sqcup_{i=1}^l L_i .$ 
Then 
\begin{align*}
\liminf_{n\to\infty}
\frac{1}{n}\log m &\left( \left\{ x\in I :
\frac{1}{n} S_n\varphi (x) > a \right\} \right)  \\
&\ge \liminf_{n\to\infty}
\frac{1}{kn}\log m
\left( \cap_{j=0}^{n-1} f^{-kn}(\sqcup_{i=1}^l L_i )\right) \\
&\ge \liminf_{n\to\infty}
\frac{1}{kn}\log \{ l^n \min_{x\in \sqcup_{i=1}^l L_i}
 |(f^k)'(x) |^{-n} m (J) \} \\
&\ge \frac{1}{k} \log l - \int \log |f' | d\mu -\varepsilon \\
&\ge h_{\mu} (f) - \int \log |f' | d\mu - 2\varepsilon 
\end{align*}
Letting $\varepsilon \to 0$
we obtain
the inequlity (\ref{ineq:les}) for the case
that $\mu$ is ergodic.
If $\mu$ is not ergodic, then
take $\varepsilon >0$ small enough so that
$\displaystyle \int \varphi d\mu >\alpha +\varepsilon $
and a linear combination
$\mu ' =\alpha_1 \mu_1 +\cdots + \alpha_p\mu_p$
of ergodic and hyperbolic  measures $\mu_1 ,\ldots , \mu_p$ 
such that
$$| h_{\mu}(f)- h_{\mu '} (f) | \le \varepsilon /4, \qquad
| \int \log |f'| d\mu - \int \log |f'| d\mu '|
 \le\varepsilon /4$$
and
$$| \int \varphi d\mu -\int \varphi d\mu ' |\le \varepsilon /2 .$$ 
Applying the argument above for each $\mu_q$, $q=1, \ldots ,p,$
we can take integers $k_q, l_q\ge 1$ with 
$(\log l_q)/k_q \ge h_{\mu_q} (f) -\varepsilon /4$  and
pairwise disjoint compact intervals $L_1^q ,\ldots , L_{l_q}^q$ with 
$L^q\subset I$
such that:
$L_i^q\subset L^q , f^{k_q}(L_i^q)=L^q$ and 
$f^{k_q}\mid_{L_i^q} :L_i^q \to L^q$ is injective on $L_i^q$
$(i=1,2,\ldots l_q)$;
$$
| \frac{1}{k_q}\log\vert (f^{k_q} )' (x) \vert 
- \int \log |f' | d\mu_q | \le\varepsilon /8
\enspace\text{ and }\enspace 
| \frac{1}{k_q} S_{k_q}\varphi (x) - \int \varphi d\mu_q | \le \epsilon /4
$$
whenever $x\in\sqcup_{i=1}^{l_q} L_i^q .$
On the other hand, 
there are integers $K_1, \ldots , K_p \ge 1$ such that
$f^{K_q} (L^q) = I ,$
since $f$ is topologically mixing. 
For any large integer $n\ge 1$ define
$r_0(n):= 0$ and 
$r_{q}(n):=r_{q-1}(n)+[n\alpha_q / k_q] k_q + K_q$
inductively on $q=1,\ldots , p,$
and set
$$B_n :=\cap_{q=1}^{p}\cap_{j=0}^{[n\alpha_q / k_q] -1}
f^{-r_{q-1}(n)-jk_q} (\sqcup_{i=1}^{l_q} L_i^q ),$$ 
where $[\cdot ]$ denotes the Gauss' symbol. 
Then since $n\ge 1$ is large, 
for any $x\in B_n$ we have 
\begin{align*}
|(f^{r_p(n)})'(x)| 
&\le 
 (\max_{y\in I} |f'(y)| )^{K_1+\cdots +K_p}
\cdot \prod_{q=1}^p |(f^{[n\alpha_q / k_q]k_q})'(f^{r_{q-1}(n)}(x))| \\
\le 
 (\max_{y\in I} &|f'(y)| )^{K_1+\cdots +K_p}
\cdot \sum_{q=1}^p \exp\left\{ [n\alpha_q / k_q]k_q 
\left(\int \log |f'| d\mu_q +\varepsilon /8 \right) \right\} \\
&\le \exp\left\{ n \left( \int \log |f'| d\mu ' + \varepsilon /4\right)
\right\} .
\end{align*}
Thus, we obtain
\begin{align*}
m(B_n)
&\ge {l_1}^{[n\alpha_1 / k_1]}\cdots {l_p}^{[n\alpha_p / k_p]}
\left(\max_{x\in B_n}  |(f^{r_p(n)})' (x) |\right)^{-1} m(I) \\
&\ge {l_1}^{[n\alpha_1 / k_1]}\cdots {l_p}^{[n\alpha_p / k_p]}
\exp \left\{ -n \left( \int \log |f'| d\mu ' + \varepsilon /4 \right)
\right\}  m(I). 
\end{align*}
Moreover, for any $x\in B_n$ we have
\begin{align*}
\Big|
S_{n}\varphi (x)&-n\int\varphi d\mu '\Big| \\
&\le \sum_{q=1}^{p}
| S_{[n\alpha_q / k_q]k_q}\varphi (f^{r_{q-1}(n)}(x))-
\left[ \frac{n\alpha_q} {k_q} \right] k_q\int\varphi d\mu_q |  \\
&+ \sum_{q=1}^p  2(k_q +K_q)\max_{y\in I}|\varphi (y) |  \\
&\le n\varepsilon /4 + n\varepsilon /4 \quad\quad
= n \varepsilon /2,
\end{align*}
and then
\begin{align*}
| \frac{1}{n} S_{n}\varphi (x)-\int\varphi d\mu | 
&\le | \frac{1}{n} S_{n}\varphi (x)-\int\varphi d\mu '| 
+ | \int\varphi d\mu '- \int\varphi d\mu |  \\
&\le \varepsilon /2 +\varepsilon /2 \quad\quad = \varepsilon .
\end{align*}
Thus,
$$\frac{1}{n} S_{n}\varphi (x) > \int\varphi d\mu - \varepsilon > a .$$
As a conclusion we have
\begin{align*}
\liminf_{n\to\infty}
\frac{1}{n}\log m \Big( \Big\{  x\in I :
\frac{1}{n} &S_n\varphi (x) > a \Big\}\Big) 
\ge \liminf_{n\to\infty}
\frac{1}{n}\log m (B_n) \\
&\ge \sum_{q=1}^p \frac{{\alpha}_q}{k_q}\log l_q 
- \left(\int \log |f'| d\mu ' + \varepsilon / 4 \right)  \\
&\ge h_{\mu '}(f)-\int \log |f'| d\mu ' - \varepsilon /2 \\
&\ge h_{\mu}(f)-\int \log |f'| d\mu  - \varepsilon .
\end{align*}
Letting $\varepsilon \to 0$
we obtain (\ref{ineq:les}), 
and hence the lower estimate of the rate function.

\noindent
{\bf The upper estimate.}
We show that for any $\varepsilon  , \eta >0$ there exists 
$\mu\in\mathcal{M}_f$ hyperbolic with $\displaystyle \int \varphi d\mu 
> a-\varepsilon $ such that
\begin{equation}
\limsup_{n\to\infty}\frac{1}{n}\log
m\Big( \Big\{ x\in I :
\frac{1}{n}S_n\varphi (x) \ge a \Big\} \Big)
\le h_{\mu} (f) -\int \log |f'|d\mu +\eta \label{ineq:ues}
\end{equation}
holds whenever the left hand side of the inequality (\ref{ineq:ues}) is not
$-\infty .$
Take a subinterval  $J\subset I$
and a return time function $R :J\to\mbox{\boldmath$N$}\cup\{ \infty \}$
as in the assumptions for $f$ stated in Section 2.  
Then setting
 $X_k :=\{ x\in J : R(x)>k \}\enspace (k=0,1,2,\ldots )$
we obtain a tower $(Z, \mathcal{A} )$ by 
$$Z:= \sqcup_{k=0}^{\infty} X_k\times \{ k\}$$
and 
$$\mathcal{A} :=  
\{ J \times \{ k \} :  J\in \mathcal{I}_k , k=0,1,2,\ldots \}$$
where $\mathcal{I}_k$ is the partition of $X_k$ 
which consists of the connected components
of both $\{ x\in J : R(x)=k+1 \}$ and $\{ x\in J : R(x)>k+1 \} .$ 
Then it follows that the tower $(Z, \mathcal{A} )$ is nonsteep  
from the assumptions for the map.
A tower map $T:Z\to Z$ defined by
$$T(x,k)  := 
\begin{cases}
(x, k+1)
 & \text{ if } R(x) > k+1  , \\
(f^{k+1} (x), 0) & \text{ if }
R(x) = k+1,
\end{cases} $$
satisfies $\pi\circ T = f\circ \pi $ on $Z$  
where $\pi (x, k)=f^k(x) .$
It also follows from the assumptions for $f$ that
the map $T:Z\to Z$ satisfies 
both of the admissibility and the bounded distortion conditions.
Moreover, if $A\in\vee_{i=0}^{n-1} T^{-i}\mathcal{A} ,$
then $\pi (A)$ is an interval with length less than 
or equal to $\varepsilon_{n} ,$ 
and it implies that the function $\psi : Z\to \mbox{\boldmath $R$}$  
defined by $\psi (x, k) =\varphi (f^k(x)) $ 
is contained in the class $\mathcal{F}$. 
Then 
since there is an integer $l\ge 1$ such that $f^lJ = I,$
taking large $n\ge 1$ we have
\begin{align*}
m\Big(\Big\{ x\in I :  &\frac{1}{n} S_n\varphi (x)
\ge a  \Big\}\Big) 
= m \Big(\Big\{ x\in f^lJ : 
 \frac{1}{n} S_n\varphi (x)
\ge a  \Bigr\}\Big) \\ 
\le &m \Big( f^l \Big(\Big\{ x\in J : 
 \frac{1}{n} S_n\varphi (x)
\ge a -\epsilon /2  \Bigr\}\Big)\Big) \\
\le &(\max_{y\in I} |(f)'(y) | )^l m\Big( \Big\{ x\in J : 
 \frac{1}{n} S_n\varphi (x)
\ge a- \epsilon /2 \Bigr\}\Big) \\
= &(\max_{y\in I} |(f)'(y) | )^l m \Big( \Big\{ x\in X_0 : 
 \frac{1}{n} S_n\psi  (x,0)
\ge a- \epsilon /2 \Bigr\}\Big) .
\end{align*}
Hence, there is $\nu\in \mathcal{M}_T$ 
with 
$\nu (\sqcup_{k=0}^{K-1} X_k\times \{ k \} )=1$
for some $K\ge 1$ such that $\displaystyle \int\psi d\nu > a-\varepsilon$
and
\begin{align*}
\limsup_{n\to\infty}\frac{1}{n}
\log 
m\Big(\Big\{ &x\in I :  \frac{1}{n} S_n\varphi (x)
\ge a  \Big\}\Big) \\
&\le 
\limsup_{n\to\infty}\frac{1}{n}\log 
 m \Big( \Big\{ x\in X_0 : 
 \frac{1}{n} S_n\psi (x,0)
\ge a- \epsilon /2 \Bigr\}\Big)  \\
&\le h_{\nu} (T)-\int \log {\rm Jac} (T)d\nu +\eta
\end{align*}
by Theorem 1. 
Then the $f$-invariant probability measure $\mu :=\nu\circ {\pi}^{-1}$ 
on $I$ is hyperbolic
because the Lyapunov exponents on the support of $\mu$ 
are not smaller than 
$(\log\lambda )/ K$ uniformly. 
Also, $\mu$ satisfies
$$\int \varphi d\mu = \int\psi d\nu > a-\varepsilon$$
and  
\begin{align*}
\limsup_{n\to\infty}\frac{1}{n}\log 
m\Big(\Big\{ x\in I :  &\frac{1}{n} S_n\varphi (x)
\ge a  \Big\}\Big) \\
&\le h_{\nu} (T)-\int \log {\rm Jac } (T) d\nu +\eta \\
&= h_{\mu} (f)-\int \log |f '| d\mu +\eta .
\end{align*}
This completes the proof of Theorem 3.


\begin{thebibliography}{99}

\bibitem{AM}{V. Ara\'ujo and M.J. Pacifico}
{\it Large deviations for non-uniformly expanding maps},
J. Stat. Phys. {\bf 125}
(2006), 411--453.


\bibitem{Ba}{V. Baladi},
{\it Positive Transfer Operators and Decay of Correlations}, 
World Scientific, Singapore, 2000.


\bibitem{BC}{M. Benedicks and L. Carleson},
{\it On iterations of $1-ax^2$ on $(-1, 1)$},
Ann. of Math. {\bf 122}
(1985), 1--25.





\bibitem{BK}{H. Bruin and G. Keller},
{\it Equilibrium states for S-unimodal maps},
Ergod. Th. Dynam. Sys. {\bf 18}
(1998), 765--789.

\bibitem{BLS}{H. Bruin, S. Luzzatto and  S. van Strien},
{\it Decay of correlations in one-dimensional dynamics},
Ann. Sci. \'Ecole Norm. Sup. {\bf 36} 
(2003), 621--646.

\bibitem{C0}{ Y.M. Chung},
{\it Shadowing property of non-invertible maps with hyperbolic measures},
Tokyo J. Math. {\bf 22}
(1999),  145--166.





\bibitem{C4}{ Y.M. Chung},
{\it Birkhoff spectra for one-dimensional maps 
with some hyperbolicity},
in preparation.



\bibitem{CE}{P. Collet and J.-P. Eckmann},
{\it Positive Lyapunov exponents and absolutely continuity
for maps of the interval},
Ergod. Th. Dynam. Sys. {\bf 3}
(1983), 13--46.

\bibitem{E}{R.S. Ellis},
{\it Entropy, Large deviations and Statistical Mechanics},
Grundlehren der Mathematischen Wissenschaften {\bf 271},
Springer-Verlag, New York, 1985.




\bibitem{K}{ A. Katok},
{\it Lyapunov exponents, entropy and periodic orbits for diffeomorphisms},
Inst. Hautes Etudes Sci. Publ. Math. {\bf 51}
(1980), 137--173.



\bibitem{KM}{ A. Katok and L. Mendoza},
{\it Dynamical systems with nonuniformly hyperbolic behavior},
supplement to "Introduction to the modern theory of dynamical 
systems" written by A. Katok and B. Hasselblatt,
Cambridge Univ. Press, Cambridge, 1995, pp 659--700.


\bibitem{Ke}{ G. Keller},
{\it Equilibrium states in ergodic theory},
London Mathematical Society Student texts {\bf 42},
Cambridge Univ. Press, Cambridge, 1998.


\bibitem{KN}{G. Keller and T. Nowicki},
{\it Spectral theory, zeta functions and 
the distribution of periodic points for 
Collet-Eckmann maps},
Commun. Math. Phys. {\bf 149}
(1992), 31--69.

\bibitem{L1}{F. Ledrappier},
{\it Some properties of absolutely continuous invariant measures
of an interval},
Ergod. Th. Dynam. Sys. {\bf 1}
(1981), 77--93. 




\bibitem{MN}{I. Melbourne and M. Nicol},
{\it Large deviations for nonuniformly hyperbolic systems},
to appear in Trans. Amer. Math. Soc. 



\bibitem{NS}{T. Nowicki and D. Sands},
{\it Non-uniform hyperbolicity and universal bounds for 
S-unimodal maps},
Invent. Math. {\bf 132}
(1998), 633--680.

\bibitem{OP}{S. Orey and S. Pelikan},
{\it Deviation of trajectory averages and the defect in Pesin's formula 
for Anosov diffeomorphisms},
Trans. Amer. Math. Soc. {\bf 315}
(1989), 741--753.



\bibitem{Pe}{Ya. Pesin},
{\it Characteristic Lyapunov exponents and smooth ergodic theory},
Russian Math. Survays {\bf 32}
(1977), 55--114.







\bibitem{PSY}{M. Pollicott, R. Sharp and M. Yuri},
{\it Large deviations for maps with indifferent fixed points},
Nonlinearity {\bf 11}
(1998), 1173--1184.


\bibitem{PM}{Y. Pomeau and P. Manneville},
{\it Intermittent transition to turbulence in dissipative dynamical systems},
Commun. Math. Phys. {\bf 74}
(1980), 189--197.


\bibitem{Ru}{D. Ruelle},
{\it An inequality for the entropy  of differentiable maps},
Bol. Soc. Brasil. Math. {\bf 9}
(1978), 83--87.

\bibitem{T0}{Y. Takahashi},
{\it Fredholm determinant of unimodal linear maps} 
Sci. Papers College Gen. Ed. Univ. Tokyo {\bf 31} (1981), 61--87. 
58F13 (28A05) 


\bibitem{T1}{Y. Takahashi}, 
{\it Entropy functional (free energy) for dynamical systems 
and their random perturbations}, 
in "Stochastic analysis (Katata/Kyoto, 1982)", 
North-Holland, Amsterdam, 1984, pp 437--467. 

\bibitem{T2}{Y. Takahashi}, 
{\it Two aspects of large deviation theory for large time}, 
in "Probabilistic methods in mathematical physics 
(Katata/Kyoto, 1985)", 
Academic Press, Boston, 1987, pp 363--384.



\bibitem{W}{P. Walters},
{\it An Introduction to Ergodic Theory},
 Graduate Texts in Mathematics {\bf 79}, Springer-Verlag, New York,
1982.


\bibitem{Y1}{L.-S. Young},
{\it Some large deviation results for dynamical systems},
Trans. Amer. Math. Soc. {\bf 318}
(1990), 525--543.

\bibitem{Y2}{L.-S. Young},
{\it Decay of correlations for certain quadratic maps},
Commun. Math. Phys. {\bf 146}
(1992), 123--138.

\bibitem{Y3}{L.-S. Young},
{\it Statistical properties of dynamical systems 
with some hyperbolicity},
Ann. of Math. {\bf 147}
(1998), 585--650.

\bibitem{Y4}{L.-S. Young},
{\it Recurrence times and rates of mixing},
Israel J. Math. {\bf 110}
(1999), 153--188.



\bibitem{Yu1}{M. Yuri},
{\it Thermodynamic formalism for certain nonhyperbolic maps},
Ergod. Th. Dynam. Sys. {\bf 19}
(1999), 1365--1378.



\end{thebibliography}
\end{document}